\documentclass[12pt,a4paper]{amsart}

\usepackage[dvips]{graphicx}

\newtheorem{theorem}{Theorem}[section]
\newtheorem{proposition}[theorem]{Proposition}
\newtheorem{lemma}[theorem]{Lemma}
\newtheorem{claim}[theorem]{Claim}

\newtheorem{corollary}[theorem]{Corollary}

\theoremstyle{definition}

\newtheorem{example}[theorem]{Example}

\theoremstyle{remark}

\numberwithin{equation}{section}

\begin{document}

\title{$P^2$-reducing and toroidal Dehn fillings}

\author[G. Jin]{Gyo Taek Jin}
\address{Department of Mathematics, Korea Advanced Institute of Science and Technology,
Taejon 305-701, Korea}
\email{trefoil@math.kaist.ac.kr}
\author[S. Lee]{Sangyop Lee}
\address{Department of Mathematics, Korea Advanced Institute of Science and Technology,
Taejon 305-701, Korea}
\email{apple@mathx.kaist.ac.kr}
\author[S. Oh]{Seungsang Oh}
\address{Department of Mathematics, Chonbuk National University, Chonju, Chonbuk 561-756, Korea}
\email{soh@math.chonbuk.ac.kr}
\author[M. Teragaito]{Masakazu Teragaito}
\address{Department of Mathematics and Mathematics Education,
Hiroshima University, Kagamiyama 1-1-1, Higashi-hiroshima 739-8524, Japan}
\email{teragai@hiroshima-u.ac.jp}



\begin{abstract}
We study the situation where we have two exceptional Dehn fillings on a
given hyperbolic $3$-manifold.
We consider two cases that one filling creates a projective plane, and
the other creates an essential torus or a Klein bottle, and
give the best possible upper bound on the distance between two fillings
for each case.
\end{abstract}
\maketitle

\section{Introduction}

Let $M$ be a compact, connected, orientable $3$-manifold with  a torus boundary
component $\partial_0M$.
A \textit{slope\/} on $\partial_0M$ is the isotopy class of an unoriented essential simple loop.
For a slope $r$, the manifold obtained from $M$ by \textit{$r$-Dehn filling\/} is
$M(r)=M\cup V_r$, where $V_r$ is a solid torus
glued to $M$ along $\partial_0M$ in such a way that $r$ bounds a meridian disk in $V_r$.
If $r$ and $s$ are two slopes on $\partial_0M$,
then $\Delta(r,s)$ denotes their minimal geometric intersection number.

We say that a $3$-manifold $M$ is \textit{hyperbolic\/} if $M$ with its boundary tori removed
admits a complete hyperbolic structure of finite volume with
totally geodesic boundary.
If $M$ has non-empty boundary,
then Thurston's geometrization theorem for Haken manifolds \cite{Th1} says that
$M$ is hyperbolic if and only if $M$ contains no essential sphere, disk, torus or annulus.

We are interested in obtaining restrictions on when a Dehn filling on a hyperbolic
$3$-manifold fails to be hyperbolic.
Such a filling is said to be \textit{exceptional}.
It is well known that if $M$ is hyperbolic, then 
there are only finitely many exceptional Dehn fillings on $\partial_0M$ \cite{Th1}, and
there are a large amount of investigations on exceptional Dehn fillings (see \cite{G3}).

In this paper, we deal with three specific exceptional Dehn fillings. 
A $3$-manifold is \textit{$P^2$-reducible\/} if it contains a projective plane, and
\textit{$P^2$-irreducible\/} otherwise.
If a closed $3$-manifold is $P^2$-reducible, then it is either the real $3$-dimensional
projective space $P^3$ or a reducible manifold with a $P^3$-summand.
A $3$-manifold is \textit{toroidal\/} if it contains an essential torus.
Clearly, if $M(r)$ is $P^2$-reducible or toroidal, then such a filling is exceptional.
Finally a hyperbolic $3$-manifold contains no Klein bottle, and
if $M(r)$ contains a Klein bottle, then the filling is also exceptional (see \cite{Lu}).

We consider two situations.
That is, one filling yields a $P^2$-reducible manifold, and the other
gives either a toroidal manifold or a manifold containing a Klein bottle.
For each case, we can find the best possible upper bound on the distance between
two exceptional fillings.

\begin{theorem}
Let $M$ be a hyperbolic $3$-manifold with a torus boundary component $\partial_0M$.
Let $\alpha$ and $\beta$ be two slopes on $\partial_0M$ such that
$M(\alpha)$ is $P^2$-reducible and $M(\beta)$ is toroidal.
Then either
\begin{itemize}
\item[(1)] $\Delta(\alpha,\beta)\le 2$; or
\item[(2)] $\Delta(\alpha,\beta)=3$ and $M(\beta)$ contains an essential torus 
which intersects the core of the attached solid torus in two points.
\end{itemize}
\end{theorem}

\begin{theorem}
Let $M$ be as in Theorem 1.1.
Let $\alpha$ and $\beta$ be two slopes on $\partial_0M$ such that
$M(\alpha)$ is $P^2$-reducible and $M(\gamma)$ contains a Klein bottle.
Then either
\begin{itemize}
\item[(1)] $\Delta(\alpha,\gamma)\le 2$; or
\item[(2)] $\Delta(\alpha,\gamma)=3$ and $M(\gamma)$ contains a Klein bottle
which intersects the core of the attached solid torus in a single point.
\end{itemize}
\end{theorem}

The examples showing that these estimates are sharp are given in the final section.

In \cite{G1}, Gordon gave an upper bound $5$ for the distance between a toroidal filling
and a lens space filling on a hyperbolic $3$-manifold with torus boundary.
Our Theorem 1.1 gives a partial improvement of this result in case where
a lens space is the real projective $3$-space $P^3$.

\begin{corollary}
Let $M$ be a hyperbolic $3$-manifold with torus boundary $\partial M$, and
let $\alpha$ and $\beta$ be slopes on $\partial M$ such that
$M(\alpha)$ is the lens space $L(2,1)\ (=P^3)$ and $M(\beta)$ is toroidal. 
Then $\Delta(\alpha,\beta)\le 3$.
\end{corollary}

The authors would like to thank Cameron Gordon for suggesting the problem.

\section{Preliminaries}

In the remainder of this paper, we assume that $M$ is a hyperbolic $3$-manifold
with a torus boundary component $\partial_0M$, and
that $\alpha, \beta, \gamma$ are slopes on $\partial_0M$ such that
$M(\alpha)$ is $P^2$-reducible, $M(\beta)$ is toroidal, and $M(\gamma)$
contains a Klein bottle.

Assume that $\Delta(\alpha,\beta), \Delta(\alpha,\gamma)\ge 3$.
Suppose that $M(\alpha), M(\beta)$ and $M(\gamma)$ contain a projective plane
$\widehat{P}$, an essential torus $\widehat{T}$ and a Klein bottle $\widehat{K}$ respectively.
Then we may assume that $\widehat{P}$ meets the attached solid torus
$V_\alpha$ in a finite collection of
meridian disks, so that $P=\widehat{P}\cap M$ is a punctured projective plane
properly embedded in $M$, each of whose boundary components has slope $\alpha$.
Furthermore, we can assume that $\widehat{P}$ is chosen so that
the number of boundary components $p$ is minimal among all projective planes in $M(\alpha)$. 
Similarly, $\widehat{T}$ and $\widehat{K}$ give rise to the surfaces
$T$ and $K$ respectively, and
the numbers of boundary components $t$ and $k$ of $T$ and $K$ are assumed to be minimal. 
Recall that $M(\alpha)$ is either the real projective $3$-space $P^3$, or
a reducible manifold with a $P^3$-summand.
Therefore we may assume that $M(\beta)$ and $M(\gamma)$ are irreducible by
\cite{BZ2,GL0}.

\begin{lemma}\label{p2}
$p\ge 2$.
\end{lemma}
 
\begin{proof}
If $p=0$, then $M$ contains a projective plane, which is impossible since $M$ is hyperbolic.
If $p=1$, then $M$ contains a M{\"o}bius band, which is also impossible.
\end{proof}

\begin{lemma}
$P$ is incompressible and boundary-incompressibe in $M$.
\end{lemma}

\begin{proof}
Assume $P$ is compressible in $M$.
Let $D$ be a compressing disk for $P$.
Note that $\partial D$ is orientation-preserving on $P$ (and hence $\widehat{P}$).
Hence $\partial D$ bounds a disk $D'$ on $\widehat{P}$.
Since $\mathrm{Int}D'$ meets $V_\alpha$, we can create a new projective plane
by replacing $D'$ with $D$, which meets $V_\alpha$ fewer than $\widehat{P}$.
This contradicts the minimality of $\widehat{P}$.
Therefore $P$ is incompressible in $M$.

Next, assume that $P$ is boundary-compressible.
Then $P$ would be compressible, or the core of $V_\alpha$ can
be isotoped into $\widehat{P}$ as an orientation-reversing loop.
But, the latter case implies that $M$ is boundary-reducible.
\end{proof}

Since $t$ is minimal, it is clear that $T$ is incompressible and boundary-incompressible in $M$.
We have $t\ge 1$, and $k\ge 1$, since $M$ cannot contain a Klein bottle.

\begin{lemma}
$K$ is incompressible and boundary-incompressible in $M$.
\end{lemma}

\begin{proof}
Suppose that $K$ is compressible in $M$.
Let $D$ be a disk in $M$
such that $D\cap K=\partial D$ and $\partial D$ does not bound a disk on $K$.
Note that $\partial D$ is orientation-preserving on $K$.

If $\partial D$ is non-separating on $\widehat{K}$, then
we get a non-separating $2$-sphere in $M(\gamma)$ by compressing $\widehat{K}$ along $D$.
This is clearly a contradiction.
If $\partial D$ bound a disk on $\widehat{K}$, then we replace the disk with $D$, and
get a new Klein bottle in $M(\gamma)$ with fewer intersections with $V_\gamma$ than
$\widehat{K}$.
This contradicts the choice of $\widehat{K}$.

Thus $\partial D$ is essential and separating on $\widehat{K}$.
Compressing $\widehat{K}$ along $D$ gives two disjoint projective planes in $M(\gamma)$.
Since $M(\gamma)$ is irreducible, this is also impossible.
Thus we have shown that $K$ is incompressible.

Next, let $E$ be a disk in $M$
such that $E\cap K=\partial E\cap K$, $\partial E=a\cup b$, where
$a\subset K$ is an essential (i.e. not boundary-parallel) arc in $K$
and $b\subset \partial_0M$.
If $a$ joins distinct components of $\partial K$, then a compressing disk for $K$ is
obtained from two parallel copies of $E$ and the disk obtained by removing a neighborhood
of $b$ from the annulus in $\partial_0M$ cobounded by those components of $\partial K$
meeting $a$. 
Hence $\partial a$ is contained in the same component $\partial_1K$, say, of $\partial K$.
If $k\ge 2$, then $b$ bounds a disk $D'$ in $\partial_0M$ together with a subarc of $\partial_1K$.
Then $E\cup D'$ gives a compressing disk for $K$ in $M$.
Therefore $k=1$. Then we can move the core of $V_\gamma$ onto an
orientation-reversing loop in $\widehat K$ by using $E$.
This implies that $M$ contains a properly embedded M{\"o}bius band,
which contradicts the fact that $M$ is hyperbolic. 
\end{proof}


From the arc components of $P\cap T$, $P\cap K$, we can construct two pairs of graphs
$(G_P^T,G_T)$ and $(G_P^K,G_K)$ in the usual way (see \cite{CGLS,G2,GL}).
We number the components of $\partial T$ as $1,2,\dots,t$ in the order in which
they appear on $\partial_0M$. Similarly number the components of $\partial K$.
But the components of $\partial P$ are numbered from $-1$ to $p-2$ unusually. 
For simplifying the notations, we use the symbol $G_P$ for both $G_P^T$ and $G_P^K$,
and $G_{TK}$ for $G_T$ or $G_K$.
For a graph $G$, the \textit{reduced graph\/} $\overline{G}$ of $G$ is defined to be the graph
obtained from  $G$ by amalgamating each family of parallel edges into a single edge.

\begin{lemma}
Neither $G_P$ nor $G_{TK}$ has trivial loops.
\end{lemma}

\begin{proof}
This follows from the fact that $P, T$ and $K$ are boundary-incompressible.
\end{proof}

We may assume that no circle component of $P\cap T$ or $P\cap K$
bounds a disk in $P, T$ or $K$,
because of the incompressibilities of these surfaces.

Although $P$ and $K$ are non-orientable, we can establish a parity rule,
which plays a crucial role in this paper.
In fact, this is a natural generalization of the usual parity rule \cite{CGLS}.

\begin{figure}[tb]
\includegraphics*[scale=1]{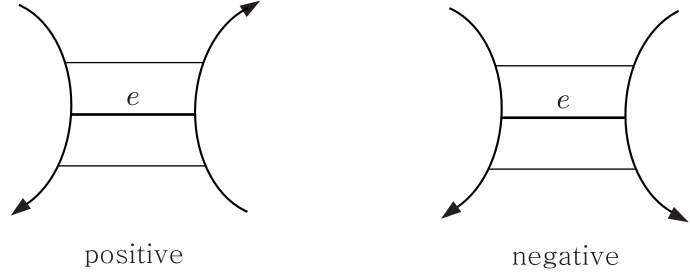}
\caption{A sign of an edge}\label{parity}
\end{figure}

First, orient all components of $\partial P$ so that they are
mutually homologous on $\partial_0M$.
Similarly for $\partial T$ and $\partial K$.
Let $e$ be an edge in $G_P$.
Since $e$ is an arc properly embedded in $P$, a regular neighborhood $D$ of $e$ in $P$
is a disk in $P$.
Then  $\partial D=a\cup b\cup c\cup d$, where $a$ and $c$ are arcs in $\partial P$ with
induced orientations from $\partial P$.
On $D$, if $a$ and $c$ have opposite directions,
then $e$ is called \textit{positive}, otherwise \textit{negative}.
See Figure \ref{parity}.
Similarly, define the sign of edges in $G_{TK}$.
Then we have the following rule.

\begin{lemma}[Parity rule]
An edge $e$ is positive \textup{(}or negative\textup{)} in $G_P$ if and
only if $e$ is negative \textup{(}positive resp.\textup{)} in $G_{TK}$.
\end{lemma}

\begin{proof}
This follows from the fact that $M$ is orientable and $\partial_0M$ is a torus.
\end{proof}

Remark that $G_T$ can have positive loops, but no negative loops.
$G_P$ and $G_K$ can have positive and negative loops.
But the key point is;

\begin{lemma}\label{loop}
\begin{itemize}
\item[(1)] At most one vertex can be a base of negative loops in $G_P$.
\item[(2)] At most two vertices can be bases of negative loops in $G_K$.
\end{itemize}
\end{lemma}

\begin{proof}
Let $e$ be a negative loop based at a vertex $x$.
Then $N(x\cup e)$ is a M{\"o}bius band.
Since $\widehat{P}$ and $\widehat{K}$ can contain at most
one M{\"o}bius band and at most two disjoint M{\"o}bius bands
respectively, the conclusions follow.
\end{proof}

Now, we have some basic properties of the graphs.
An edge is called an \textit{$x$-edge\/} if it has label $x$ at its endpoint, and
a \textit{level edge\/} if both endpoints have the same label.
For example, $G_{KT}$ can have a positive level edge, which corresponds to
a negative loop in $G_P$ by the parity rule.
Thus at most one label of $G_{TK}$ can be a label of positive level edges by Lemma \ref{loop}(1).
Therefore, we adopt the convention that $0$ is the label of positive level edges
in $G_{TK}$.

Let $G=G_P$ or $G_{TK}$.
A cycle in $G$ is a \textit{Scharlemann cycle\/}
if it bounds a disk face, and the edges in the cycle are all positive and
have the same label pair $\{i,i+1\}$ at their two endpoints,
called the \textit{label pair\/} of the Scharlemann cycle.
(In this case, the label set of $G$ must have at least two elements.)
In particular, a Scharlemann cycle of length two is called an \textit{$S$-cycle\/} in short.

When $p\ge 3$, a \textit{generalized Scharlemann cycle\/} in $G_{TK}$ is defined to be
either a Scharlemann cycle or
a cycle of positive edges whose labels are in the set $\{-1,0,1\}$ at both ends and
which bounds a disk face in $G_{TK}-\{\mathrm{level\ edges}\}$.
In addition, if $p=3$, then the label $0$ must appear around each vertex in the disk face
bounded by the generalized Scharlemann cycle in $G_{TK}-\{\mathrm{level\ edges}\}$.
For example, Figure 2(b) is not a generalized Scharlemann cycle.
Refer to \cite{DM} for more details.

A generalized Scharlemann cycle of length two is called
a \textit{generalized $S$-cycle}.
Then a generalized $S$-cycle, not an $S$-cycle, in $G_{TK}$ is
a triple $\{e_{-1},e_0,e_1\}$ of mutually
parallel positive edges where $e_{-1}$ and $e_1$ have the same label pair $\{-1,1\}$,
and $e_0$ is a level edge with label $0$.
A generalized $S$-cycle, not an $S$-cycle,
can be defined in $G_P$ (precisely, $G_P^K$),
but then the label set is $\{i-1,i,i+1\}$ for some $i$.

\begin{figure}[tb]
\includegraphics*[scale=1]{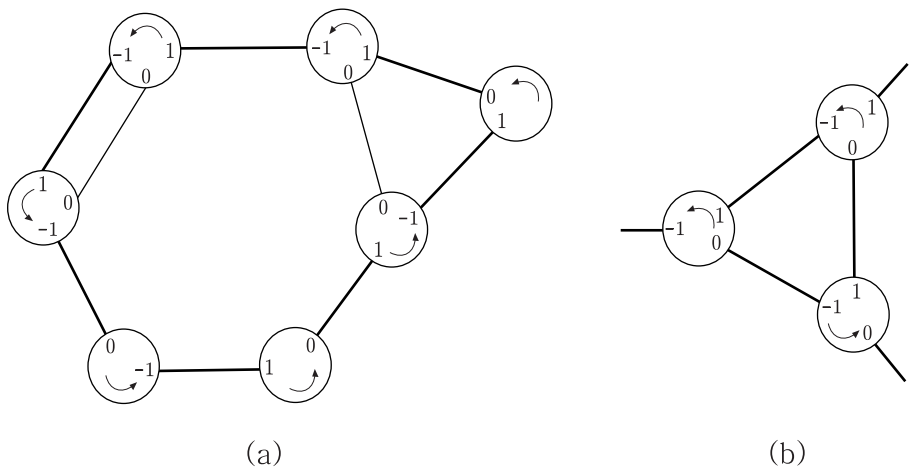}
\caption{A generalized Scharlemann cycle}\label{fig1}
\end{figure}

\begin{lemma}\label{SS}
 $G_{TK}$ contains no generalized Scharlemann cycles.
\end{lemma}

\begin{proof}
%
This is \cite[Theorem 1.1]{DM}.
\end{proof}

\begin{lemma}\label{positive}
\begin{itemize}
\item[(1)] If $p$ is even, $G_{TK}$ has at most $p/2$ mutually parallel positive edges.
Furthermore, if there is a family of $p/2$ mutually parallel positive edges, then
either the first, say, edge of the family is level, or
the set of labels of one end of the family is disjoint from that of the other end.
\item[(2)] If $p$ is odd, $G_{TK}$ has at most $(p+1)/2$ mutually parallel positive edges.
Furthermore, if there is a family of $(p+1)/2$ mutually parallel positive edges, then
the first, say, edge of the family is level.
\end{itemize}
\end{lemma}

\begin{proof}
These follow from Lemmas \ref{loop}(1) and \ref{SS}.
\end{proof}

\begin{lemma}\label{p-parallel}
Let $p\ge 2$. Then $G_{TK}$ cannot contain $p$ mutually parallel edges.
\end{lemma}

\begin{proof}
Let $A_1,A_2,\dots,A_p$ be a family of  $p$ mutually parallel edges
in $G_{TK}$ labelled successively.
By Lemma \ref{positive}, all $A_i$'s are negative, and then
they make orientation-preserving cycles in $G_P$.
Note that an orientation-preserving loop in a projective plane is contractible.
Therefore, we can choose an innermost cycle among them.
Then the construction in \cite[Section 5]{GLi} implies that $M$ is cabled, a contradiction.
\end{proof}

\section{Main argument}

In this section, we use $|G_{TK}|$ to denote the number of vertices of $G_{TK}$.

For a label $x$ of $G_{TK}$, let $\Gamma_x$ be the subgraph of $G_{TK}$ consisting of all the vertices
and positive $x$-edges of $G_{TK}$.
A disk face $D$ of $\Gamma_x$ is called an \textit{$x$-face} of $G_{TK}$.

Let $G_P^+$ denote the subgraph of $G_P$ consisting of all the vertices of $G_P$ and
the positive edges of $G_P$.
Note that $G_P^+$ has a disk support in $\widehat{P}$, that is, there is a disk in $\widehat{P}$
which contains $G_P^+$ in its interior, since any orientation-preserving loop
in a projective plane is contractible.

Let $\Lambda$ be a subgraph of $G_P^+$ with a disk support $D$. 
A vertex of $\Lambda$ is a \textit{boundary vertex\/} if
there is an arc connecting it to $\partial D$ whose interior is disjoint from
$\Lambda$,
and an \textit{interior vertex\/} otherwise.

A \textit{generalized web\/} $\Lambda_P$ is
a connected subgraph of $G_P^+$ satisfying that
\begin{itemize}
\item[(i)] at most one boundary vertex $y$ of $\Lambda_P$ is a cut vertex of $G_P^+$; 
\item[(ii)] each vertex of $\Lambda_P$, except $y$ if it exists,
            has degree at least $(\Delta-1)|G_{TK}|$ in $\Lambda_P$;
\item[(iii)] there is a disk $D$ in $\widehat{P}$ such that $D\cap G_P^+=\Lambda_P$.
\end{itemize}

The vertex $y$ as in (i), if it exists, is called an \textit{exceptional vertex} of
$\Lambda_P$.

\begin{proposition}\label{gweb}
One of the followings holds\textup{;}
\begin{itemize}
\item[(1)] $G_{TK}$ contains an $x$-face for some $x\ne 0$\textup{;}
\item[(2)] $G_P$ contains a generalized web.
\end{itemize}
\end{proposition}

\begin{proof}
We distinguish two cases.

\bigskip
\noindent
\textbf{Case 1:} \enspace Suppose that there is a vertex $x (\ne 0)$ of $G_P$
such that more than $|G_{TK}|$ negative edges are incident to $x$.
Remark that such negative edges are not loops because $x\ne 0$.
This implies that there exist more than $|G_{TK}|$ positive $x$-edges in $G_{TK}$ by the parity rule.
Thus $\Gamma_x$ has a larger number of edges than that of vertices.
An Euler characteristic calculation gives that $\Gamma_x$ contains a disk face.

\bigskip
\noindent
\textbf{Case 2:} \enspace As the negation of Case 1, suppose that
each vertex $x (\ne 0)$ of $G_P$ has at least $(\Delta-1)|G_{TK}|$ positive edge endpoints.

Let $\Lambda$ be an extremal component of $G_P^+$.
That is, $\Lambda$ is a component of $G_P^+$ having a disk support $D$ such that
$D\cap G_P^+=\Lambda$.

First, assume that $\Lambda$ is a single vertex.
Then it must be vertex $0$.
Thus only negative edges are incident to vertex $0$ in $G_P$.
If there is no negative loop at vertex $0$, then $\Lambda$ can be the only extremal component of $G_P^+$.
But then, $G_{TK}$ contains a $0$-face $E$ as in Case 1, and $G_{TK}$ has no positive level edges.
Then $E$ contains a Scharlemann cycle by \cite[Proposition 5.1]{HM}, which
is impossible by Lemma \ref{SS}.
If a negative loop is incident there, we can choose another extremal component of $G_P^+$,
which has more than one vertex.
Thus we can assume that $\Lambda$ is not a single vertex.

Choose a block $\Lambda_P$ of $\Lambda$ with at most one cut vertex.
Then $\Lambda_P$ is clearly a generalized web.
\end{proof}

\section{$x$-face of $G_{TK}$}

In this section, we treat the case (1) of Proposition \ref{gweb}.

\begin{theorem}\label{gene}
Let $p\ge 3$.
If $G_{TK}$ contains a non-zero $x$-face $D$, then it contains
a generalized Scharlemann cycle in $D$.
\end{theorem}

\begin{proof}
There is a possibility that $\partial D$ is not a circle.
That is, $\partial D$ may contain a double edge, and also
more than two edges of $\partial D$ may be incident to a vertex on $\partial D$. 
Since we will find a generalized Scharlemann cycle within $D$,
we can cut formally the graph $G_{TK}\cap D$ along double edges of $\partial D$
and at vertices to which more than two edges of $\partial D$ are incident so that
$\partial D$ is deformed into a circle.
(See also \cite[Section 5]{HM}.)
Thus we may assume that $\partial D$ is a circle.

If $D$ is a bigon, then the conclusion is obvious.
Therefore, we assume that $D$ has at least three sides.

Suppose that $D$ has a diagonal edge $e$, which has two distinct labels
$a$ and $b$ at its endpoints as in Figure \ref{split}(a).
Since $D$ is an $x$-face, $a\ne x$ and $b\ne x$.
Without loss of generality, we may assume that the labels
appear in counterclockwise order around the boundary of each vertex, 
and that $a>b$.

\begin{figure}[tb]
\includegraphics*[scale=1]{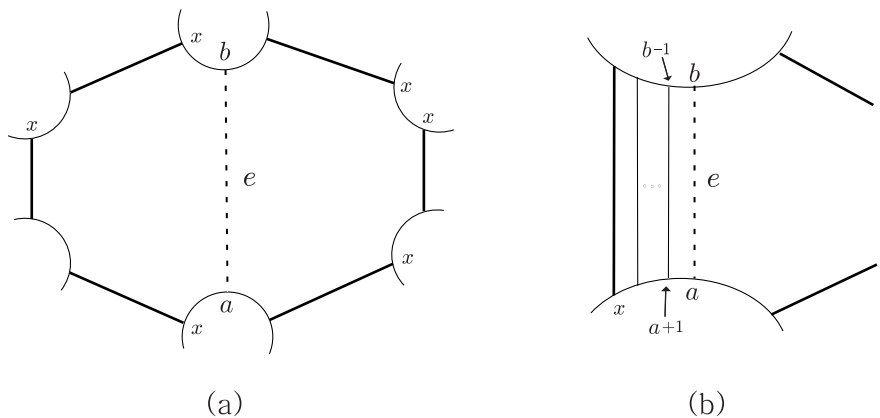}
\caption{Split along a diagonal edge}\label{split}
\end{figure}

Formally, construct a new $x$-face $D'$ as follows.
The edge $e$ divides $D$ into two disks $D_1$ and $D_2$.
We can assume that $D_1$ lies on the right side of $e$ when $e$ is oriented from
the endpoint with label $a$.
If three labels $b,a,x$ appear in this order around the corners of $\partial D$,
then discard $D_2$, and 
insert additional edges to the left of $e$, and parallel to $e$, until we
first reach label $x$ at one or both ends of an additional edge.
See Figure \ref{split}(b).
The new $x$-face $D'$ is the union of $D_1$ and some additional bigons.
If three label $b,x,a$ appear in this order, then discard $D_1$, and
insert additional edges to the right of $e$ as above.
Then $D'$ is the union of $D_2$ and some additional bigons.
Remark that there is no generalized Scharlemann cycle among additional edges and $e$.

Repeat the above process for every diagonal edge which is not level,
then get a new $x$-face $E$ and a graph $\Gamma$ in $E$.
All diagonal edges of $\Gamma$ are level, and all boundary edges are $x$-edges,
especially label $x$ can appear on both ends of a boundary edge.
Such boundary edges are called level $x$-edges to distinguish them from level $0$-edges.
See Figure \ref{gamma}.

\begin{figure}[tb]
\includegraphics*[scale=1]{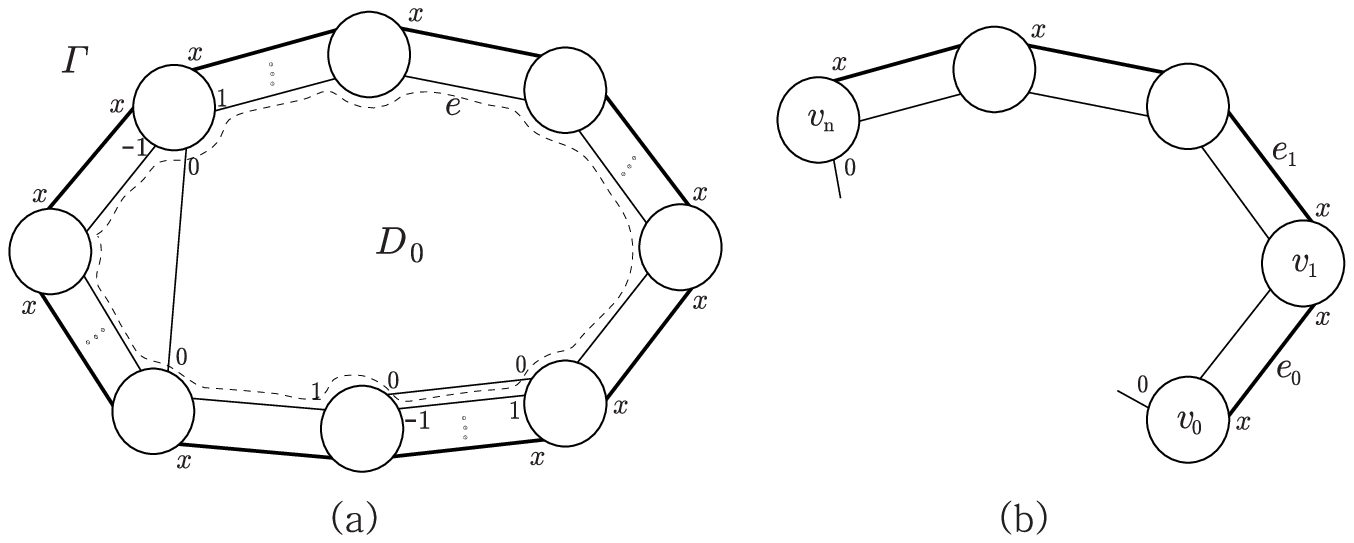}
\caption{$\Gamma$}\label{gamma}
\end{figure}

From now on, we assume that there is no Scharlemann cycle in $\Gamma$.

\begin{claim}
$\Gamma$ contains a level $0$-edge.
\end{claim}

\begin{proof}
Suppose not.
Then $\Gamma$ contains a Scharlemann cycle by \cite[Lemma 5.2]{HM}.
\end{proof}

Let $D_0$ be the disk face of
$\Gamma-\{ \mbox{level 0-edges} \subset \mathrm{Int}E\}$,
which is not a bigon.
We will show that for any edge of $\Gamma$ on $\partial D_0$, it has
labels only on $\{-1,0,1\}$ at both ends.

There are consecutive vertices $v_0,v_1,\dots,v_n$ on $\partial D_0$ such that
$v_i$ is not a base of a level $0$-edge for $1\le i\le n-1$, but
$v_0$ and $v_n$ are base of level $0$-edges.
Possibly, $n=1$.  See Figure \ref{gamma}(b).
Let $e_i$ be the $x$-edge connecting $v_i$ with $v_{i+1}$ for $i=0,1,\dots,n-1$.
Let $r(e_i)$ ($l(e_i)$, resp.)  be the label of the end of $e_i$ at $v_i$ ($v_{i+1}$, resp.)
and $F_i$ the family of mutually parallel edges containing $e_i$.  
The number of edges in $F_i$ is denoted by $|F_i|$.
Also, let $f_i$ be the edge of $F_i$ lying on $\partial D_0$.
We define $l(f_i)$ and $r(f_i)$ similarly as above.

\begin{claim}\label{F}
\begin{itemize}
\item[(1)] If $p$ is even, then $|F_i|\le p/2$.
Furthermore, if $|F_i|=p/2$, then 
either the $x$-edge $e_i$ is level, or its label on the other end is $x-1$ or $x+1$.
\item[(2)] If $p$ is odd, then $|F_i|\le (p+1)/2$.
Furthermore, if $|F_i|=(p+1)/2$, then
$e_i$ is a level $x$-edge.
\end{itemize}
\end{claim}

\begin{proof}
These follow from Lemma \ref{positive} and the definition of $F_i$.
\end{proof}

\begin{claim}\label{one}
If $x=\pm 1$, or $|F_j|=1$ for some $j$, then
each $f_i$ has labels on $\{-1,0,1\}$.
\end{claim}

\begin{proof}
Suppose $x=1$.
Then $r(e_0)\ne 1$, since $F_0$ cannot have $p-1$ edges, except the case $p=3$,
by Claim \ref{F}.
Thus $l(e_0)=1$ and $|F_0|\le 2$, since there cannot be a
level edge or an $S$-cycle.
When $p=3$, we have the same conclusion.
For, $|F_0|=1$ or $2$, and if two then $l(e_0)=r(e_0)=1$ by Claim \ref{F}.

The same argument runs until we got that $l(e_i)=1$ and $|F_i|\le 2$ for all $i$.
Then it is easy to see that $f_i$ has labels on $\{-1,0,1\}$ at both ends (see Figure \ref{pm1}(a)).
In particular, when $p=3$, the corner on $\partial D_0$ at $v_i$ contains
label $0$.
Similarly for the case where $x=-1$.

\begin{figure}[tb]
\includegraphics*[scale=1]{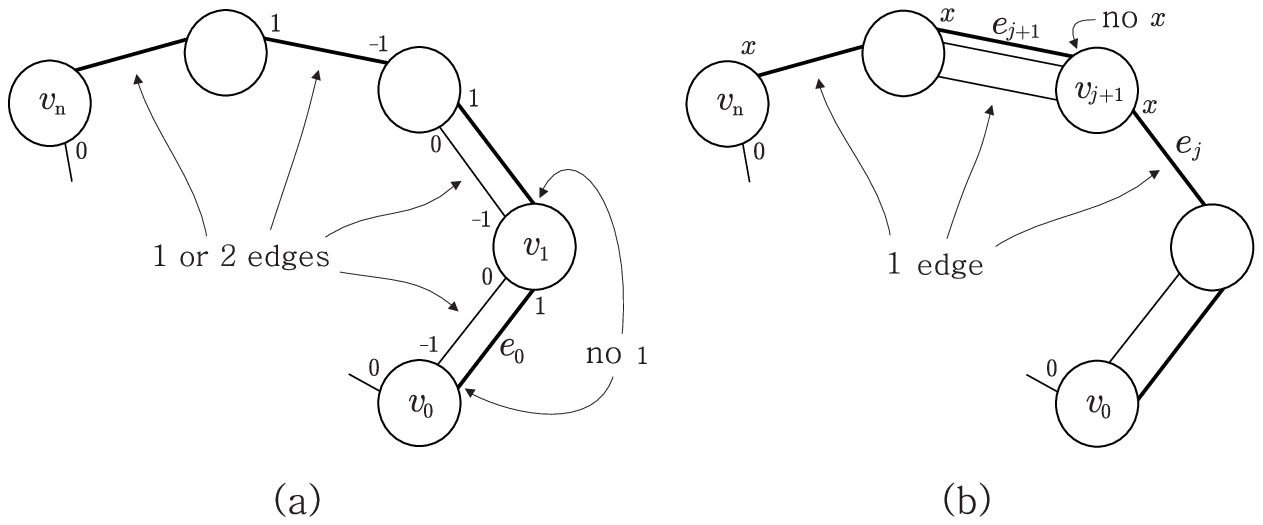}
\caption{}\label{pm1}
\end{figure}

Suppose that $|F_j|=1$ as in Figure \ref{pm1}(b).
Assume $l(e_j)=x$.
By Claim \ref{F}, $r(e_{j+1})\ne x$, and so $l(e_{j+1})=x$.
If $|F_{j+1}|>1$, then
it contains either a level edge or an $S$-cycle.
Thus $|F_{j+1}|=1$.
Hence $|F_{n-1}|=1$ and $l(e_{n-1})=x$.
This implies that $x=1$.
If $r(e_j)=x$, then a similar argument shows that $|F_0|=1$ and $r(e_0)=x$,
which implies $x=-1$.
\end{proof}

When $p=3$, then $x=\pm 1$.
Then we get the desired result by Claim \ref{one}.
Therefore, we assume $p>3$.
Furthermore, we can assume that $x\ne\pm 1$ and any $|F_i|>1$ by Claim \ref{one}.

\begin{claim}
If $r(e_i)=x$ for all $i$, then each $f_i$ has labels on $\{-1,0,1\}$.
\end{claim}

\begin{proof}
Note that $r(f_0)=-1$ and $l(f_{n-1})=1$.
Consider $F_0$.
Since $F_0$ contains neither level edge nor $S$-cycle,
$l(f_0)\ge 0$. 
This implies $|F_0|\le |F_1|$.
Then $l(f_1)\ge l(f_0)$.
Thus we have
$0\le l(f_0)\le l(f_1)\le \dots \le l(f_{n-1})=1$.
Hence $l(f_i)=0$ or $1$ for all $i$.
The result immediately follows from this observation.
\end{proof}

Thus we can assume that $r(e_i)=x, 0\le i\le m-1$ and $r(e_m)\ne x$ for some $m$.
(Possibly, $m=0$.)  Then $l(e_m)=x$.
See Figure \ref{red}.
Also, $|F_m|\le p/2$ by Claim \ref{F}.
Then $r(e_{m+1})\ne x$.
Thus we have $r(e_i)\ne x$ and 
$l(e_i)=x$ for $m\le i\le n-1$, and $|F_{n-1}|\le p/2$.

\begin{figure}[tb]
\includegraphics*[scale=1]{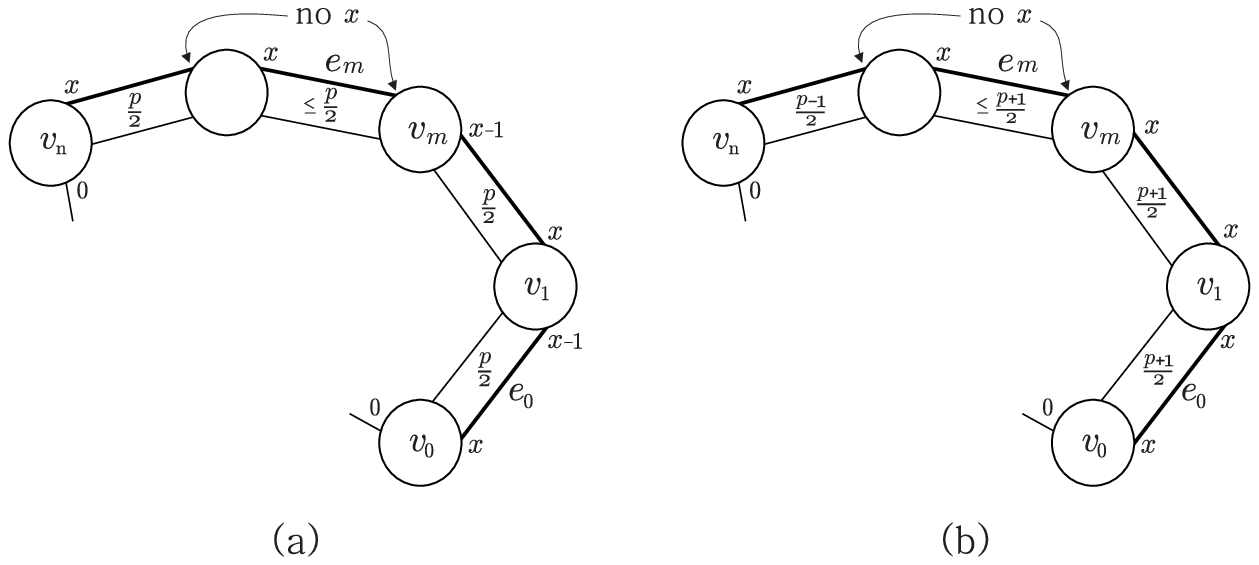}
\caption{}\label{red}
\end{figure}

\bigskip
\noindent
\textbf{Case 1:}\enspace \textit{$p$ is even.}

Then $|F_0|=|F_{n-1}|=p/2$, since
$|F_0|+|F_{n-1}|=p$.  See Figure \ref{red}(a).
Thus $l(e_0)=x-1$ or $x$ by Claim \ref{F}.
For, if $l(e_0)=x+1$ then $F(e_0)$ contains an $S$-cycle.
Indeed, $l(e_0)=x-1$, since $r(e_1)=x$.
Hence $|F_1|=p/2$, and $l(e_1)=x-1$ or $x$ again.
If $r(e_2)=x$, then $l(e_1)=x-1$ as above.
Thus we can conclude that $|F_i|=p/2$ $(0\le i\le m-1)$,
and that $l(e_i)=x-1$, and $l(f_i)=0$, $r(f_i)=-1$ for $0\le i\le m-2$.
Also $l(e_{m-1})=x$ or $x-1$, and hence $l(f_{m-1})=0$ or $1$, and $r(f_{m-1})=-1$. 
Since $|F_{n-1}|=p/2$ and $l(f_{n-1})=1$,
we see $r(e_{n-1})=x+1$ and so $r(f_{n-1})=0$.
Thus $l(f_i)=1$ for $m\le i\le n-1$, and
$r(f_i)=0$ for $m+1\le i\le n-1$ and $r(f_m)=-1$ or $0$.

\bigskip
\noindent
\textbf{Case 2:}\enspace \textit{$p$ is odd.}

Then $|F_0|=(p+1)/2$, and $|F_{n-1}|=(p-1)/2$ similarly for the case where $p$ is even.
See Figure \ref{red}(b).
Also, $l(e_0)=x$ by Claim \ref{F}(2).
Then we see that $|F_i|=(p+1)/2$ and $l(e_i)=x$ for $0\le i\le m-1$.
Thus $r(f_i)=-1$ and $l(f_i)=0$ for $0\le i\le m-1$.
Since $|F_{n-1}|=(p-1)/2$ and $l(f_{n-1})=1$,
we have $r(e_{n-1})=x+1$ and $r(f_{n-1})=-1$.
Then $l(f_{n-2})=0$ and therefore $|F_{n-2}|=(p+1)/2$.
By Claim \ref{F}(2), the only possibility is $m=n-1$.

Thus we have shown that $\partial D_0$ gives a generalized Scharlemann cycle.
This completes the proof of Theorem \ref{gene}.
\end{proof}

\begin{theorem}\label{web}
Let $p\ge 3$.  Then 
$G_{TK}$ cannot contain a non-zero $x$-face.
\end{theorem}

\begin{proof}
This follows immediately from Lemma \ref{SS} and Theorem \ref{gene}. 
\end{proof}


\section{The case $G_{TK}=G_K$}

\begin{lemma}\label{propK}
\begin{itemize}
\item[(1)] $G_P$ cannot contain a Scharlemann cycle.
\item[(2)] $G_P$ cannot contain a generalized $S$-cycle.
\item[(3)] $G_P$ has at most two level edges with different labels.
\item[(4)] Suppose $k\ge 3$.  Then $G_P$ has at most $k/2+1$
\textup{(}$(k+1)/2$, resp.\textup{)}
mutually parallel positive edges if $k$ is even \textup{(}odd, resp.\textup{)}.
Furthermore, if $G_P$ has a family of $k/2+1$ \textup{(}$(k+1)/2$, resp.\textup{)}
mutually parallel positive edges, then
the first and last edge \textup{(}the first or last edge, resp.\textup{)}
of the family are level, when $k$ is even \textup{(}odd, resp.\textup{)}.
\end{itemize}
\end{lemma}

\begin{proof}
If $G_P$ contains a Scharlemann cycle, then we can get a new Klein bottle
in $M(\gamma)$ which meets $V_\gamma$ fewer than $\widehat{K}$ (\cite[Theorem 6.4]{GL}).
Theorem 6.6 in \cite{GL} implies (2).
Here, we need the fact that the distance between two Dehn fillings creating
projective planes is at most one \cite{Te}. 
(3) is a restatement of Lemma \ref{loop}(2).  
For (4), if not, such a family contains a generalized $S$-cycle.
\end{proof}

\begin{lemma}\label{pge3}
If $k\ge 2$, then $p \ge 3$.
\end{lemma}

\begin{proof}
Assume for contradiction that $p=2$ because of Lemma \ref{p2}.
Let $u$ and $v$ be the vertices of $G_P$, where $u$ can be
a base of negative loops.
Since $u$ and $v$ have the same degree $\Delta k$, if $u$ has a loop,
then so does $v$.
Then there would be a trivial loop.  Thus we can see that $G_P$ has no loops.
Then $G_P$ consists of at most two families of mutually parallel edges;
one is a family of positive edges, and the other is that of negative edges.

If $k\ge 3$, then $G_P$ contains more than $k$ mutually parallel negative edges by Lemma  
\ref{propK}(4).
Then an easy Euler characteristic calculation shows that $G_K$ contains
a $0$-face and no level edges.
As in Case 2 of the proof of Proposition \ref{gweb}, this gives a contradiction.

If $k=2$, then $G_P$ has at most two positive edges.
Otherwise, there would be two edges which are parallel in both $G_P$ and $G_K$.
But this implies that $M$ is cabled by \cite[Lemma 2.1]{G2}.
Thus $G_P$ contains at least four negative edges.
Similarly, $G_K$ contains a $0$-face and no level edges, a contradiction.
\end{proof}

\begin{theorem}\label{nowebk}
If $k\ge 3$, then 
$G_P$ cannot contain a generalized web. 
\end{theorem}

\begin{proof}
Assume for contradiction that $G_P$ contains a generalized web $\Lambda_P$, 
possibly with an exceptional vertex $y$ among boundary vertices of $\Lambda_P$.
Let $D$ denote a disk support of $\Lambda_P$.
Lemma \ref{propK}(3) guarantees the existence of a label $x$ such that
$G_P$ contains no positive level $x$-edges.

Consider $\Lambda_P^x$ consisting of all vertices and $x$-edges of $\Lambda_P$.
Since every boundary vertex of $\Lambda_P$, except $y$, has degree at least $(\Delta - 1)k$,
it has at least two edges attached with label $x$.
We remark that $\Lambda_P^x$ may be disconnected.
Choose an innermost component $G$ of $\Lambda_P^x$ (in $D$), and
let $H$ be its block with at most one cut vertex of $G$.

Let $v,e$ and $f$ be the numbers of vertices, edges, and disk faces
of $H$, respectively. 
(We view $H$ as the graph in a disk.)
Also let $v_i, v_{\partial}$ and $v_c$ be the numbers
of interior vertices, boundary vertices of $H$ and a cut vertex of $G$ in $H$, respectively.
Hence $v=v_i+v_{\partial}$ and $v_c=0$ or $1$. 

Since $H$ has neither a level $x$-edge nor a generalized $S$-cycle, 
each face of $H$ is a disk with at least $3$ sides. 
Thus we have $3f + v_{\partial} \le 2e$. 
Combined with $1=v-e+f$ because it has only disk faces, 
we get $e \le 3v_i + 2v_{\partial} - 3$.
On the other hand we have $2(v_{\partial} - v_c) + 3v_i \le e$
because each boundary vertex of $H$, except a cut vertex of $G$,
has at least two edges attached with label $x$.  
These two inequalities give us that $3 \le 2v_c$, a contradiction.
\end{proof}

\section{The case $k=2$}

Consider the case where $k=2$ and $\Delta\ge 3$ when $G_{TK}=G_K$.
Recall that $G_P$ contains a generalized web $\Lambda_P$
by Proposition \ref{gweb}, Theorem \ref{web} and Lemma \ref{pge3}.
We remark that at most two positive edges can be parallel in $G_P$.
For, if there are three mutually parallel positive edges, then 
$G_P$ contains either an $S$-cycle or a pair of positive level edges with the same label.
The former is impossible by Lemma \ref{propK}.
In the latter case, such two edges are also parallel in $G_K$.
But this implies that $M$ is cabled \cite[Lemma 2.1]{G2}.
Also, if there are two parallel positive edges, then both edges must be level.

A positive level edge with label $1$ (or $2$) in $G_P$
is called a \textit{$1$-edge\/} (\textit{$2$-edge}, resp.)
and a positive non-level edge is called a \textit{mixed edge}.

\begin{lemma}\label{pair}
$\Lambda_P$ contains a pair of parallel edges.
In particular, these edges are a $1$-edge and a $2$-edge.
\end{lemma}

\begin{proof}
Any interior vertex of $\Lambda_P$ has degree at least six, since $\Delta k\ge 6$.
Also, any boundary vertex of $\Lambda_P$, except an exceptional one,
has degree at least four.
Therefore $\Lambda_P$ contains a pair of parallel edges by
\cite[Lemmas 2.3, 3.2]{Wu}.
Then these edges are both level as above.
\end{proof}

\begin{lemma} \label{3level}
Let $x$ be a vertex of $\Lambda_P$.
\begin{itemize}
\item[(1)] There are no three $i$-edges at $x$ for each $i=1, 2$.
\item[(2)] There are no three \textup{(}positive\textup{)} level edges connecting $x$
to mutually distinct vertices, not $x$.
\end{itemize}
\end{lemma}

\begin{proof}
(1) Suppose that three positive $1$-edges, say, are incident to $x$.
Since all negative loops based at vertex $1$ in $G_K$ are mutually parallel,
there are at least $p+1$ negative loops at vertex $1$.
This is impossible by Lemma \ref{p-parallel}.

(2) Let $v_1$ and $v_2$ be the vertices of $G_K$.
By Lemma \ref{pair}, $\Lambda_P$ contains a pair of parallel (positive) edges $e_1, e_2$.
Here, $e_i$ is an $i$-edge for $i=1,2$, and $e_i$ connects a vertex $a$ to another $b$ (possibly, $a=b$).
Then $e_i$ gives a negative loop based at $v_i$ with the label pair $\{a,b\}$.
Remark that all negative loops at $v_i$ are mutually parallel.
If $\Lambda_P$ has another $i$-edge $e$ connecting $x$ to $y$, then
$e_i$ and $e$ are parallel, and so $y=x\pm (b-a) \pmod{p}$.
This means that any $i$-edge at $x$ goes to either vertex $x+b-a$ or $x-b+a$.
\end{proof}

Let $n_i$ denote the number of negative loops at vertex $i$ in $G_K$.
Without loss of generality, we can assume that $n_1\ge n_2$.

\begin{lemma} \label{mark}
Let $x$ be a vertex of $\Lambda_P$.
If a level edge $e$ and a pair of parallel edges are incident to $x$ successively in $\Lambda_P$,
then $e$ is a $1$-edge.
Furthermore, there are no consecutive pairs of parallel edges at $x$ in $\Lambda_P$.
\end{lemma}

\begin{proof}
Let $e_1, e_2$ be these pair edges, where
$e_i$ is an $i$-edge.
Assume for contradiction that
$e$ is the consecutive $2$-edge next to $e_1$.
As in the proof of \ref{3level}(2), we may assume that
$e$ connects $x$ to $x+r$ and $e_1$ and $e_2$ connect $x$ to $x-r$ for some $r$.
Let $\alpha, \beta, \alpha'$ and $\beta'$ be the endpoints of these edges.
Let $\overrightarrow{\alpha\beta}$ be a properly oriented arc from $\alpha$ to $\beta$ 
along $\partial P$, and $\overrightarrow{\alpha'\beta'}$ similarly.  
See Figure \ref{key}(b).

\begin{figure}[tb]
\includegraphics*[scale=1]{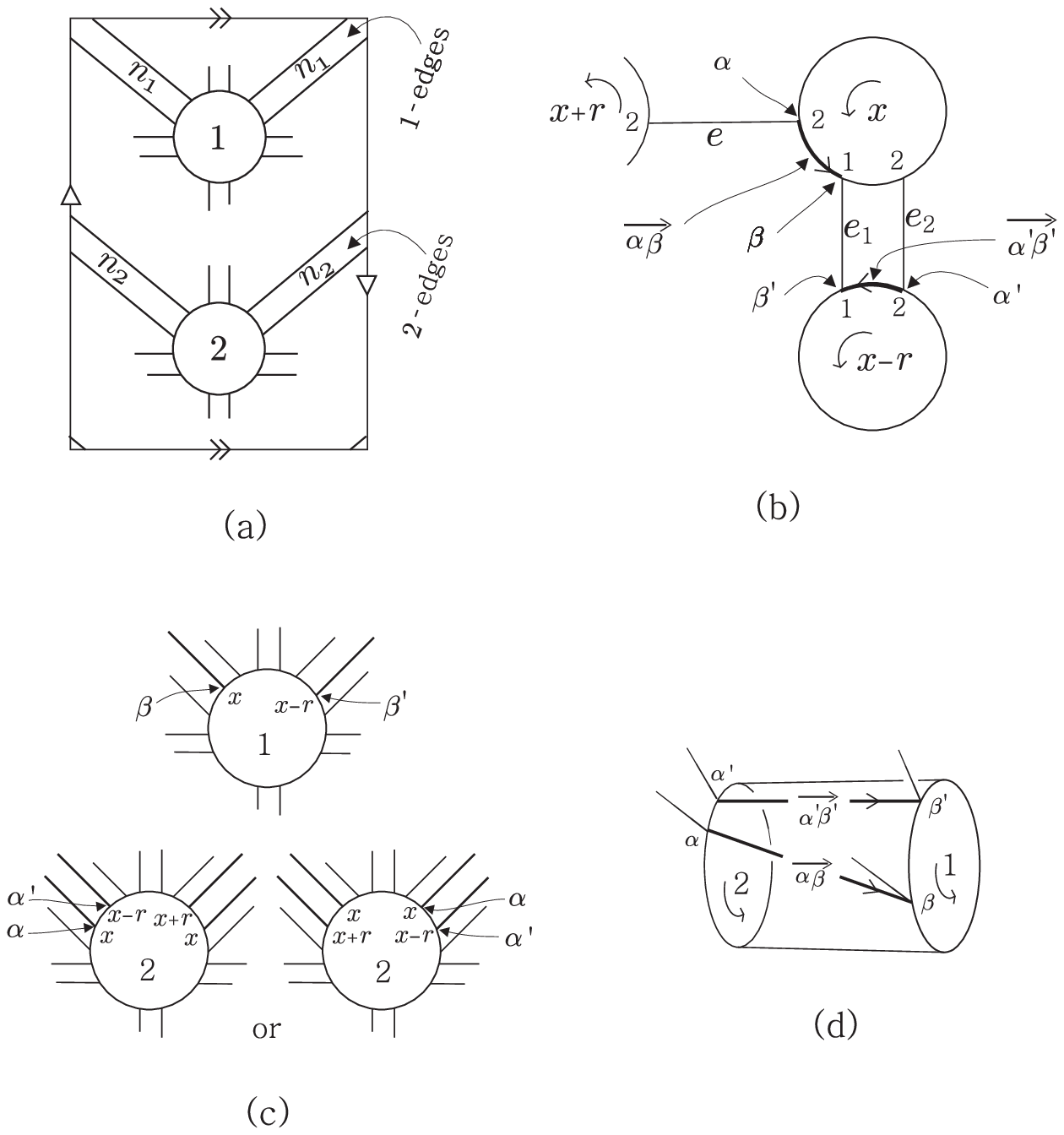}
\caption{}\label{key}
\end{figure}

In $G_K$, $e_1$ is a negative loop at vertex $1$, and
its endpoint $\beta$ has label $x$, and $\beta'$ has $x-r$.
Also, $e$ and $e_2$ are contained in the family of mutually parallel
negative loops at vertex $2$ in $G_K$, and
$\alpha$ and $\alpha'$ appear at the same end of this family.
Otherwise, the family contains more than $p$ edges, contradicting Lemma \ref{p-parallel}.

Let $\partial_iK$ be the boundary component of $K$ with label $i$ ($i=1,2$).
We may assume that $\partial_2K$ is oriented
counterclockwise as in Figure \ref{key}(c).
Then the subarc $\alpha'\alpha$, with the induced orientation from $\partial_2K$,
contains at most $n_2$ edge endpoints.


Consider the annulus part of $\partial_0M$ between $\partial_1K$ and $\partial_2K$,
containing $\overrightarrow{\alpha\beta}$ and $\overrightarrow{\alpha'\beta'}$.
Then we see that the two oriented subarcs $\alpha'\alpha$ and $\beta'\beta$ contain
the same number of points in $\partial P\cap \partial K$.  See Figure \ref{key}(d).
But, the subarc $\beta'\beta$ contains at least $n_1+1$ edge endpoints.

The second conclusion immediately follows from the first one.
\end{proof}

\begin{lemma}\label{mix}
$\Lambda_P$ has the following properties.
\begin{itemize}
\item[(1)] Every disk face is bounded by mixed edges and a non-zero even number of level edges.
\item[(2)] There is no disk bounded by a cycle consisting of mixed edges,
except at most one edge, and containing no vertex in its interior.
\end{itemize}
\end{lemma}

\begin{proof}
(1) If a disk face is bounded by only mixed edges, then the boundary cycle is
a Scharlemann cycle, which is impossible by Lemma \ref{propK}(1).
If the boundary cycle contains an odd number of level edges,
we cannot put the labels correctly.

(2) Suppose that there is such a disk $D$. 
Since its interior contains no vertex of $\Lambda_P$,
there is an outermost disk face of $\Lambda_P$ whose boundary edges
consist of mixed edges and at most one diagonal edge (or exceptional one), contradicting (1).
\end{proof}

Suppose that $\Lambda_P$ has no interior vertex.
Then, $\Lambda_P$ has a boundary vertex $x$ of degree
two in $\overline{\Lambda}_P$ by [24, Lemma 3.2].
Then there are two consecutive pairs of parallel edges incident to $x$ in $\Lambda_P$.  
This contradicts Lemma \ref{mark}.
Thus $\Lambda_P$ has an interior vertex. 
Lemma \ref{3level}(1) implies that each interior vertex has at least two mixed edges.

Consider a maximal path consisting of mixed edges starting at an interior vertex.
By Lemma \ref{3level}(2), both ends will reach some boundary vertices.
Thus we have a consecutive sequence of mixed edges between two boundary vertices.
Then we can choose an outermost disk $D$ so that its interior contains
no vertex and no mixed edges, and its boundary consists of a consecutive sequence of 
mixed edges and a consecutive sequence of boundary edges.
Furthermore, we can assume that $D$ does not contain
an exceptional vertex of $\Lambda_P$.

\begin{figure}[tb]
\includegraphics*[scale=1]{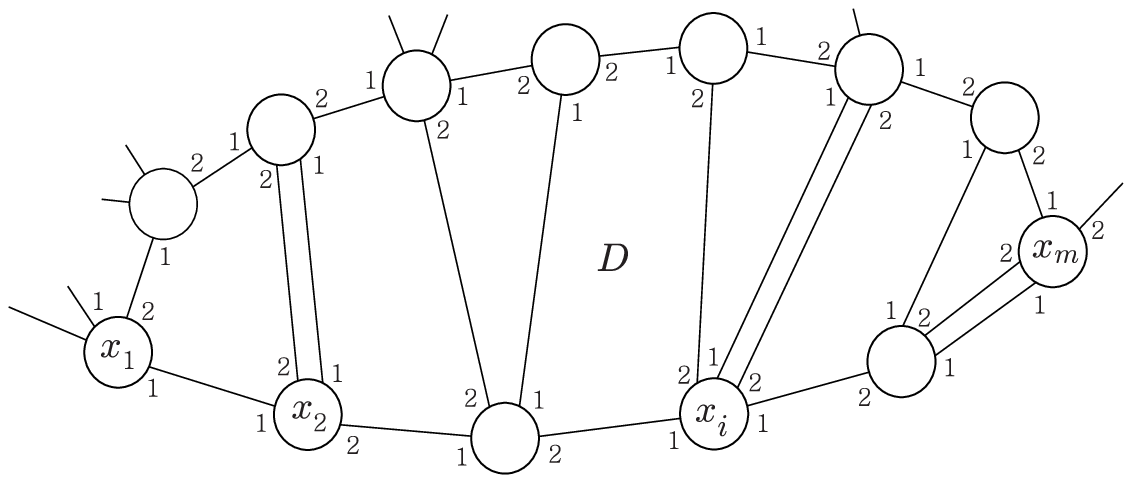}
\caption{}\label{last}
\end{figure}

In $D$, every diagonal edge of $\overline{\Lambda}_P$ connects an interior vertex and
a boundary vertex of $\Lambda_P$. 
For, if there is a diagonal edge connecting two interior vertices, 
then it contradicts Lemma \ref{mix}(2). 
Also, if there is a diagonal edge connecting two boundary vertices, 
then between these two boundary vertices there must be a vertex 
which has degree two in $\overline{\Lambda}_P$, contradicting Lemma \ref{mark}.
Recall that boundary vertices has degree at least four in $\Lambda_P$. 
See Figure \ref{last}.

Let $x_1, x_2, \cdots , x_m$ be the boundary vertices from left to right. 
There are at least three boundary vertices by Lemma \ref{mix}(2) and each vertex 
$x_i, i=2,\cdots,m-1$, has diagonal edges (so level edges) incident to by Lemma \ref{mark}.
Thus all boundary edges between $x_2$ and $x_{m-1}$ are mixed edges, 
and the other boundary edges are level edges by Lemma \ref{mix}(1).
Furthermore vertices $x_2$ and $x_{m-1}$ have a pair of parallel edges and a consecutive level edge $e$. 
By Lemma \ref{mark}, $e$ is a $1$-edge. But there must
be another vertex $x_i$ between $x_2$ and $x_{m-1}$ 
which has three consecutive level edges, so that the level edge not of the pair is a $2$-edge.
This contradicts Lemma \ref{mark}.


\begin{theorem}
If $k\ge 1$, then $\Delta\le 2$.
\end{theorem}

\begin{proof}
By Lemma \ref{pge3}, $p\ge3$.
If $k\ge 3$, then Proposition \ref{gweb}, Theorems \ref{web} and \ref{nowebk} give a contradiction.
We have just shown that the case $k=2$ is
impossible.
\end{proof}

\section{The case $k=1$}

\begin{theorem}
If $k=1$, then $\Delta\le 3$.
\end{theorem}

\begin{proof}
Assume that $\Delta\ge 4$.
The reduced graph $\overline{G}_K$
is a subgraph of the graphs shown in Figure \ref{reduced}.

\begin{figure}[tb]
\includegraphics*[scale=1]{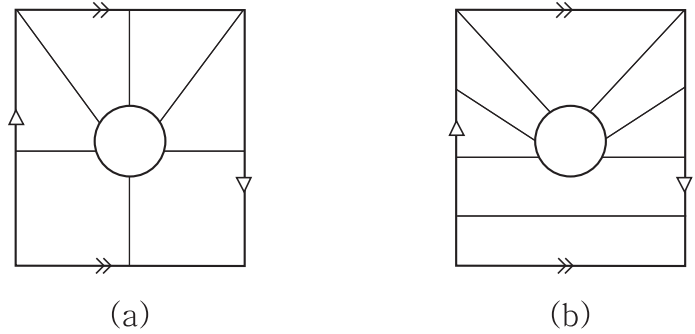}
\caption{$\overline{G}_K$}\label{reduced}
\end{figure}

If $p$ is even, there are at most $p/2$ positive edges by Lemma \ref{positive}, and
at most $p-1$ mutually parallel negative edges by Lemma \ref{p-parallel}. 
Thus $\Delta p\le  p+4(p-1)=5p-4$, and then $\Delta=4$ and $p\ge 4$.
Similarly, if $p$ is odd, we have that $\Delta=4$ and $p\ge 3$.
It follows that $G_K$ contains at least two positive edges.

First, assume that $\overline{G}_K$ is a subgraph of Figure \ref{reduced}(a).  
This family of mutually
parallel positive edges in $G_K$ contains a generalized $S$-cycle, contradicting Lemma \ref{SS}.
Next, assume that $\overline{G}_K$ is a subgraph of Figure \ref{reduced}(b).   
Since each family of mutually parallel negative edges contains at most $p-1$ edges,
we can assume that the labels are as in Figure \ref{b}.
Then the family of mutually parallel positive edges contain an $S$-cycle with label pair $\{p-2,-1\}$.
\end{proof}

\begin{figure}[tb]
\includegraphics*[scale=1]{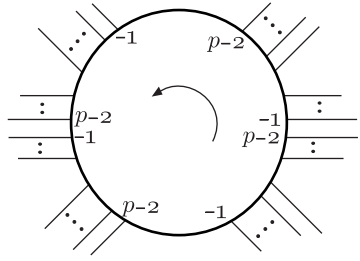}
\caption{Labels in $G_K$}\label{b}
\end{figure}

Thus we have proved Theorem 1.2.

\section{The case $G_{TK}=G_T$}

An \textit{extended $S$-cycle} in $G_P$ is the quadruple $\{e_1,e_2,e_3,e_4\}$ of
mutually parallel positive edges in succession and $\{e_2, e_3\}$ forms an
$S$-cycle.

\begin{lemma}\label{SSt}
Let $t\ge 3$.
\begin{itemize}
\item[(1)] $G_P$ cannot contain an extended $S$-cycle.
\item[(2)] $G_P$ has at most $3$ labels which are labels of $S$-cycles.
\item[(3)] $G_P$ has at most $t/2+1$ mutually parallel positive edges. 
\end{itemize}
\end{lemma}

\begin{proof}
(1) is shown in \cite[Lemma 2.10]{BZ1}.

(2) If there are four such labels, there are two $S$-cycles with disjoint
label pairs.  Then $M(\beta)$ contains a Klein bottle as in the proof of \cite[Lemma 3.10]{GL}.
(Recall that $M(\beta)$ is irreducible.) 
By Theorem 1.2, $M(\beta)$ must contain a Klein bottle $F$ which meets the core $k_\beta$ of
the attached solid torus $V_\beta$ in a single point, since $\Delta\ge 3$.
Then the torus $\partial N(F)$ is essential in $M(\beta)$.
Otherwise $M(\beta)$ would be either a lens space containing a Klein bottle
or a prism manifold, both of which are not toroidal.
This essential torus meets $k_\beta$ in two points.  This contradicts the minimality of $t$.

(3) follows from (1) and (2), since a family of $t/2+2$ mutually parallel positive edges
contains either an extended $S$-cycle or two $S$-cycles with disjoint label pairs.
\end{proof}

\begin{lemma}\label{pge3t}
If $t\ge 3$, then $p\ge 3$.
\end{lemma}

\begin{proof}
This follows from the same argument as in the proof of Lemma \ref{pge3}.
We use Lemma \ref{SSt}(3) instead of Lemma \ref{propK}(4).
\end{proof}

\begin{theorem}\label{nowebt}
If $t\ge 3$, then
$G_P$ cannot contain a generalized web.
\end{theorem}

\begin{proof}
Assume for contradiction that $G_P$ contains a generalized web $\Lambda_P$, 
possibly with an exceptional vertex $y$ of $G_P^+$.
Let $D$ denote a disk support of $\Lambda_P$.

Let $x$ be a label of $G_P$.
Suppose that $G_P$ has no $S$-cycle with a label $x$.
Consider $\Lambda_P^x$, consisting of all vertices and $x$-edges of $\Lambda_P$
as in the proof of Theorem \ref{nowebk}.
Choose an innermost component $G$ of $\Lambda_P^x$ (in $D$), and
let $H$ be its block with at most one cut vertex of $G$.
Since $H$ has neither $S$-cycle with label $x$ nor extended $S$-cycle, 
each face of $H$ is a disk with at least $3$ sides. 
Then the same calculation as in the proof of Theorem \ref{nowebk} gives a contradiction.
Therefore $G_P$ has $t$ labels which are labels of S-cycles.
Note that the existence of an S-cycle in $G_P$ guarantees
that $t$ is even \cite[Lemma 2.2]{BZ1} and so $t\ge 4$.
This contradicts Lemma \ref{SSt} (2).
\end{proof}

\begin{theorem}
If $t\ge 3$, then $\Delta\le 2$.
\end{theorem}

\begin{proof}
By Lemma \ref{pge3t}, $p\ge3$.
Then Proposition \ref{gweb}, Theorems \ref{web} and \ref{nowebt} give a contradiction.
\end{proof}

\section{The case $t=2$}  

\begin{theorem}
If $t=2$, then $\Delta\le 3$.
\end{theorem}

\begin{proof}
The reduced graph $\overline{G}_T$ is a subgraph of the graph
shown in Figure \ref{5-1} \cite[Lemma 5.2]{G2}.

\begin{figure}[tb]
\includegraphics*[scale=1]{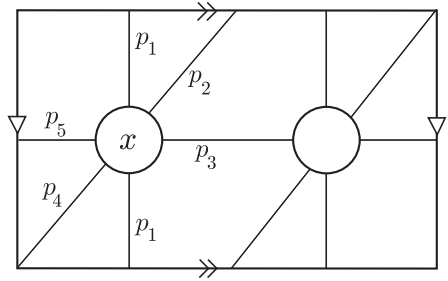}
\caption{$\overline{G}_T$}\label{5-1}
\end{figure}

Here, $p_i\ge 0$ denotes the number of edges in each of the families of
mutually parallel edges.
Then $\Delta p=2p_1+p_2+p_3+p_4+p_5$.
We have that $p_1\le (p+1)/2$ and $p_i\le p-1$ for $i=2,3,4,5$ by Lemmas \ref{positive} and \ref{p-parallel}.
Thus, $\Delta p\le (p+1)+4(p-1)=5p-3$.
This inequality implies that $\Delta=4$ and all $p_i$'s are non-zero.

\begin{claim}
All non-loop edges of $G_T$ are negative.
\end{claim}

\begin{proof}
Assume not.
Let $x$ and $y$ be the vertices of $G_T$.
If $p$ is even, then $4p=\Delta p \le 6 p/2=3p$, which is absurd.
If $p$ is odd, then $4p=\Delta p\le 6 (p+1)/2=3p+3$.
Thus $p=3$.  But then, each of six families at $x$, say, consists of
two edges, and we have an $S$-cycle in the family of loops at vertex $x$.
\end{proof}

Without loss of generality, we can assume that $p_1+p_2+p_3\ge 2p$.
Let $r\equiv p_1+p_2+p_3+1 \pmod{p}$.
Since $2p<p_1+p_2+p_3+1\le (p+1)/2+2(p-1)+1<3p$, we see $p_1+p_2+p_3+1=r+2p$. 
Then $p_1=r+2p-(p_2+p_3+1)\ge r+2p-2(p-1)-1=r+1$.
Hence $1\le r\le p_1-1$.
Thus the loop family around $x$, say, contains a generalized $S$-cycle.
\end{proof}

\section{The case $t=1$}  

\begin{theorem}
If $t=1$, then $\Delta \le 1$.
\end{theorem}

\begin{proof}
Assume that $\Delta\ge 2$.
Since $t=1$, $G_T$ can have only positive edges, 
and so $G_P$ has only negative edges.
Note that $G_T$ has at most three families of mutually 
parallel edges by [8, Lemma 5.1].
Let $A$, $B$ and $C$ be such three families, 
and $|A|$, $|B|$ and $|C|$ be the number of edges of each family.
Without loss of generality, we can assume that $A$ has at least one edge.

First assume that $\Delta$ is even.
If $|A| > 1$, then $A$ contains a generalized S-cycle by the symmetry
of labels around the vertex. Thus $|A| =1$.
Similarly we have $|B|\le 1$, $|C| \le 1$,
and therefore $\Delta p \le 6$.
Since $\Delta$ is even and $p\ge 2$,
we have $\Delta=2$ and $p=2,3$.
Then $G_T$ has two or three level edges with different labels,
contradicting Lemma 2.6(1).

Assume now that $\Delta$ is odd. 
Since $\Delta p /2$ is integral, $p$ must be even.
Thus $|A| \le p/2$ by Lemma \ref{positive}, and similarly for $B, C$.
Then $\Delta p/2=|A|+|B|+|C|\le 3p/2$.
Hence we have $\Delta=3$ and $|A|=|B|=|C|=p/2$.
This implies that $G_T$ contains two Scharlemann cycles of length three, 
which is a contradiction.
\end{proof}

We have thus completed the proof of Theorem 1.1.

\section{Examples}

In this section we give examples showing that our estimates are sharp.

\begin{example}
\cite[Theorem 4.2]{EW} shows that there is a hyperbolic
manifold $M$, admitting two Dehn fillings $M(\alpha)$ and $M(\beta)$ such that
$M(\alpha)=L(3,1)\sharp L(2,1)$ and $M(\beta)$ is toroidal, and $\Delta(\alpha,\beta)=3$.
In fact, $M(\beta)$ contains an essential torus which bounds $Q(2,-4)$ and $Q(2,2)$,
where $Q(r,s)$ is a Seifert fibered manifold with orbifold a disk with two cone points of
index $r$ and $s$.
Also, $Q(2,2)$ contains a Klein bottle, and so does $M(\beta)$. 
It is not hard to see that $M(\beta)$
contains an incompressible torus hitting the attached solid torus $V_\beta$ twice, and
a Klein bottle hitting $V_\beta$ once from the construction using tangles.
\end{example}

\bibliographystyle{amsplain}

\end{document}